\documentclass[11pt,a4paper]{article}
\usepackage[top=25mm,bottom=35mm,left=20mm,right=20mm,truedimen]{geometry}

\usepackage{hyperref}

\usepackage{amsmath,newtxtext,newtxmath,bm}

\usepackage{amsthm}
\usepackage{mathtools}
\mathtoolsset{showonlyrefs=true}
\usepackage{graphicx}
\newtheorem{jthrm}{Theorem}
\newtheorem{jlmm}[jthrm]{Lemma}
\newtheorem{jdef}[jthrm]{Definition}
\newtheorem{jpro}[jthrm]{Proposition}
\newtheorem{jcor}[jthrm]{Corollary}
\newtheorem*{jprf}{Proof}
\usepackage{commands}
\newcommand{\mqty}[1]{\begin{bmatrix}#1\end{bmatrix}}

\title{Relation between the T-congruence Sylvester equation and the generalized Sylvester equation}
\author{
  Yuki Satake\thanks{Graduate School of Engineering, Nagoya University, Furo-cho, Chikusa-ku, Nagoya 464-8603, Japan, Email: \texttt{\{y-satake,m-oozawa,sogabe,kemmochi,zhang\}@na.nuap.nagoya-u.ac.jp}},~
  Masaya Oozawa$^\ast$,~
  Tomohiro Sogabe$^\ast$\\
  Yuto Miyatake\thanks{Cybermedia Center, Osaka University, 1-32 Machikaneyama, Toyonaka, Osaka 560-0043, Japan,~ Email: \url{miyatake@cas.cmc.osaka-u.ac.jp}},~
  Tomoya Kemmochi$^\ast$,~
  Shao-Liang Zhang$^\ast$
}
\date{}

\begin{document}
\maketitle
\begin{abstract}
  The T-congruence Sylvester equation is the matrix equation $AX+X^{\mathrm{T}}B=C$, where $A\in\mathbb{R}^{m\times n}$, $B\in\mathbb{R}^{n\times m}$, and $C\in\mathbb{R}^{m\times m}$ are given, and $X\in\mathbb{R}^{n\times m}$ is to be determined.
  Recently, Oozawa et al.\ discovered a transformation that the matrix equation is equivalent to one of the well-studied matrix equations (the Lyapunov equation); 
  however, the condition of the transformation seems to be too limited because matrices $A$ and $B$  are assumed to be square matrices ($m=n$).     
  In this paper, two transformations are provided for rectangular matrices $A$ and $B$. One of them is an extension of the result of Oozawa et al. for the case $m\ge n$, and 
  the other is a novel transformation for the case $m\le n$.
\end{abstract}
\section{Introduction}
      We consider the T-congruence Sylvester equation of the form
    \begin{align} \label{T-Sylvester}
      AX+X^{\mathrm{T}}B=C,
    \end{align}
    where $A\in\mathbb{R}^{m\times n}$, $B\in\mathbb{R}^{n\times m}$, and $C\in\mathbb{R}^{m\times m}$ are given, and $X\in\mathbb{R}^{n\times m}$ is to be determined.
    The matrix equation~\eqref{T-Sylvester} appears in the industrial and  applied mathematical field: palindromic eigenvalue problems
    for studying the vibration of rail tracks under the excitation arising from high-speed trains~\cite{byers2006, kressner2009} and also in an abstract mathematical field: the tangent spaces of the Stiefel manifold~\cite{absil2009}.
    
\par
In the past decade, the numerical solvers for solving equation~\eqref{T-Sylvester} in view of direct methods and iterative methods have been well developed, (see e.g.,~\cite{chiang2012, teran2011, hajarian2014}); however, the theoretical features of the matrix equation itself  such as the uniqueness of the solution have not been explored extensively.  
In 2018, Ter\'an et al.  provided an important theoretical contribution regarding the necessary and sufficient condition for the existence of exactly one solution of the matrix equation~\eqref{T-Sylvester} in~\cite{teran2018}, whose theoretical approach is based on transforming the matrix equation into large linear systems and then analyzing the linear systems.
  \par
Though irrelevant to the above, Oozawa et al., in 2018, developed a novel transformation that can be used to transform the matrix equation \eqref{T-Sylvester} into a more well-researched matrix equation, i.e., $X^{\mathrm{T}}$ does not appear in its transformed matrix equation, as shown in the following theorem:

    \begin{jthrm}{\upshape (\cite[Theorem 3.1]{oozawa2018})}\label{moth}
      Let $A, B, C\in\mathbb{R}^{n\times n}$ and assume that $A$ is nonsingular.
      Moreover, let $S:=B^{\mathrm{T}}A^{-1}$ and suppose $\Lambda(S)$ is reciprocal free.
      Then the T-congruence Sylvester equation (1) is equivalent to the Lyapunov equation
      \begin{align}\label{lyap_oozawa}
          \tilde{X}-S\tilde{X}S^{\mathrm{T}}=Q,
      \end{align}
      where $\tilde{X}:=AX$ and $Q:=C-(SC)^{\mathrm{T}}$.
    \end{jthrm}

The significance of Theorem~\ref{moth} is that finding mathematical features and solvers of the matrix equation \eqref{T-Sylvester} reduces to only finding them of the Lyapunov equation whose features are studied extensively in control theory~\cite{gajic2008} and whose solvers have been well developed (see e.g.,~\cite{datta2004}).  
Furthermore, very recently, Miyajima utilized Theorem~\ref{moth}  to design a fast verified algorithm for the solution of the matrix equation~\eqref{T-Sylvester} in~\cite{miyajima2018}.
From the viewpoint of the complexities of the numerical error verification, utilizing Theorem~\ref{moth}  successfully reduces the complexities of $O(n^6)$ to $O(n^3)$.
\par
   Although Theorem~\ref{moth} is of interest, the class of matrices that satisfies the required condition is limited.
In fact,  Theorem~\ref{moth} does not hold if matrices $A$ and $B$ are rectangular.  
This motivates us to find transformations for the case where matrices $A$ and $B$ are rectangular. We will see in the following sections that 
a transformation for the case $m\ge n$ can be regarded as an extension of Theorem~\ref{moth}, and  a transformation for the case $m\le n$ is  a novel 
transformation in that a similar mathematical approach as used in Theorem~\ref{moth} cannot be employed mainly due to under-determined linear systems underlying the matrix equation.
\par
  The rest of this paper is organized as follows.
  Section 2 describes some mathematical preliminaries on the Kronecker product.
  Section 3 presents the main results that under some conditions the T-congruence Sylvester equation is  equivalent to a generalized Sylvester equation.
  The conclusion is provided in Section 4.
\par
    Throughout this paper, $\Lambda(S)$ denotes the spectrum of $S$, i.e., the set of eigenvalues of $S\in\mathbb{R}^{m\times m}$.
    Let $I_m$ denote the $m\times m$ identity matrix and $A^\dagger$ the Moore--Penrose inverse of $A$.
    \section{Preliminaries}
      In this section,  some mathematical preliminaries are collected on the Kronecker product and the vectorization of a matrix, referred to as {\it the  vec operator}.
      Given $A=[a_{ij}]\in\mathbb{R}^{m\times n}$ and $B\in\mathbb{R}^{p\times q}$, the Kronecker product $\otimes$ is defined by
      \begin{equation}\label{kronecker_product}
        A\otimes B := \mqty{a_{11}B & a_{12}B & \cdots & a_{1n}B \\a_{21}B & a_{22}B & \cdots & a_{2n}B \\\vdots & \vdots & \ddots & \vdots \\a_{m1}B & a_{m2}B & \cdots & a_{mn}B}\in\mathbb{R}^{mp\times nq}.
      \end{equation}
      In addition, for $C\in\mathbb{R}^{n\times l}$ and $D\in\mathbb{R}^{q\times r}$, it is known that
      \begin{equation}\label{kronecker_thrm}
        (A\otimes B)(C\otimes D) = (AC)\otimes (BD).
      \end{equation}
      \par Any eigenvalue of $A\otimes B$ is a product of an eigenvalue of $A$ and an eigenvalue of $B$ as shown below.
      \begin{jlmm}{\upshape (\cite[p.412]{lancaster1985})} \label{kronecker_eigenvalue}
        Let $\lambda_1,\ldots,\lambda_m$ be the eigenvalues of $A\in\mathbb{R}^{m\times m}$ and $\mu_1,\ldots,\mu_n$ be the eigenvalues of $B\in\mathbb{R}^{n\times n}$. 
        Then the eigenvalues of $A\otimes B$ can be written as
        \begin{equation}
          \lambda_i \mu_j \quad (1\leq i\leq m, 1\leq j\leq n).
        \end{equation}
      \end{jlmm}
      For $A=[\bm{a}_1, \bm{a}_2, \ldots, \bm{a}_n]\in\mathbb{R}^{m\times n}$, the vec operator, $\mathrm{vec}: \mathbb R^{m\times n}\to\mathbb R^{mn}$ is defined by
      \begin{equation}
        \mathrm{vec}(A):=\mqty{\bm{a}_1\\ \bm{a}_2\\ \vdots\\ \bm{a}_n}\in\mathbb{R}^{mn},
      \end{equation}
      and the inverse vec operator, $\mathrm{vec}^{-1}:\mathbb{R}^{mn}\to\bRmn$ is defined by
      \begin{equation}
        \mathrm{vec}^{-1}(\mathrm{vec}(A))=A.  
      \end{equation}
      For the vec operator, the following fact holds:
      \begin{jlmm}{\upshape (\cite[Theorem 10]{zhang2013})}\label{vec_ABC}
        Let $A\in\mathbb{R}^{m\times n}$, $B\in\mathbb{R}^{n\times p}$, and $C\in\mathbb{R}^{p\times q}$. Then, it follows that
        \begin{equation}
          \mathrm{vec}(ABC)=(C^{\mathrm{T}}\otimes A)\mathrm{vec}(B).
        \end{equation}
      \end{jlmm}
      We collect some properties on permutation matrices related to the vec operator.
      \begin{jlmm}{\upshape (\cite{zhang2013})}\label{permu_matrix}
        Let $\bm{e}_{in}$ denote an $n$-dimensional column vector which has 1 in the $i$th position and 0 for otherwise, i.e.,
        \begin{equation}
          \bm{e}_{in}:=[0,0,\ldots,0,1,0,\ldots,0]^{\mathrm{T}}\in\mathbb{R}^{n}.
        \end{equation}
        Let $A\in\mathbb{R}^{m\times n}$ and $B\in\mathbb{R}^{p\times q}$.
        Then, for the permutation matrix
        \begin{equation}
          P_{mn}:=\mqty{I_m \otimes \bm{e}_{1n}^{\mathrm{T}} \\
              I_m \otimes \bm{e}_{2n}^{\mathrm{T}} \\
              \vdots \\
              I_m \otimes \bm{e}_{nn}^{\mathrm{T}}
            }\in\mathbb{R}^{mn\times mn},
        \end{equation}
        the following properties hold:
        \begin{align}
          &P_{mn}^{\mathrm{T}} = P_{nm},\label{permu_1} \\
          &P_{mn}^{\mathrm{T}}P_{mn} =P _{mn}P_{mn}^{\mathrm{T}}=I_{mn},\label{permu_2} \\ 
          &\mathrm{vec}(A) = P_{mn}\mathrm{vec}(A^{\mathrm{T}}), \label{permu_3} \\
          &P_{mp}(A\otimes B)P_{nq}^{\mathrm{T}} = B\otimes A. \label{permu_4}
        \end{align}
      \end{jlmm}
      It is known that T-congruence Sylvester equation~\eqref{T-Sylvester} can be transformed into a linear system whose coefficient matrix has 
      the Kronecker product structure with the permutation matrix multiplication ($P_{mm}$). To see this, 
      applying the vec operator to~\eqref{T-Sylvester} with Lemma~\ref{vec_ABC} yields
      \begin{align}\label{vec_T-Sylvester}
        (I_m\otimes A)\mathrm{vec}(X)+(B^{\mathrm{T}}\otimes I_m)\mathrm{vec}(X^{\mathrm{T}})=\mathrm{vec}(C).
      \end{align}
      It follows  from~\eqref{permu_2},~\eqref{permu_3}, and~\eqref{permu_4} of Lemma~\ref{permu_matrix} that  the second term of the left-hand side of equation~\eqref{vec_T-Sylvester} is
      \begin{align}
        (B^{\mathrm{T}}\otimes I_m)\mathrm{vec}(X^{\mathrm{T}})&=(B^{\mathrm{T}}\otimes I_m)P_{mn}\mathrm{vec}(X)\\
        &=P_{mm}(I_m\otimes B^{\mathrm{T}})P_{mn}^{\mathrm{T}}P_{mn}\mathrm{vec}(X)\\
        &=P_{mm}(I_m\otimes B^{\mathrm{T}})\mathrm{vec}(X).
      \end{align}
      Thus, T-congruence Sylvester equation~\eqref{T-Sylvester} can be transformed into the following linear system:
      \begin{align}
        \{I_m\otimes A+P_{mm}(I_m\otimes B^{\mathrm{T}})\}\bm{x}=\bm{c},\label{linear_eq}
      \end{align}
      where $\bm{x}:=\mathrm{vec}(X)$ and $\bm{c}:=\mathrm{vec}(C)$.
    \section{Main results}
      
      In this section, a transformation of the T-congruence Sylvester equation into a generalized Sylvester equation is shown. A transformation for the case  $m\geq n$ is provided in Section 3.1, which is an extension of Theorem~\ref{moth},  and  a novel transformation for the case $m\leq n$ is provided in Section 3.2.
      Before stating the main results, the term {\it reciprocal free} is described.
      The following definition admits 0 and $\infty$ as elements of $\lambda_1, \ldots , \lambda_n$.
      \begin{jdef}\label{reciprocal_free}
        {\upshape (reciprocal free~\cite{byers2006,kressner2009})} A set $\{\lambda_1, \ldots, \lambda_n\}\subset\mathbb{C}\cup\{\infty\}$ is said to be reciprocal free if $\lambda_i\neq \frac{1}{\lambda_j}$ for any $1\leq i, j\leq n$.
      \end{jdef}
      \subsection{The case $m\geq n$}
      In this subsection, a transformation from the T-congruence Sylvester equation~\eqref{T-Sylvester} to the generalized Sylvester equation for the case $m\geq n$ is described.
      \begin{jthrm}\label{overdetermined}
        Let $m\geq n$, $A\in\mathbb{R}^{m\times n}$, $B\in\mathbb{R}^{n\times m}$, and $C\in\mathbb{R}^{m\times m}$. 
        Assume that there exists a matrix $S\in\mathbb{R}^{m\times m}$ such that $B^{\mathrm{T}}=SA$ and $\Lambda(S)$ is reciprocal free.
        Then, the T-congruence Sylvester equation~\eqref{T-Sylvester} is equivalent to the generalized Sylvester equation
        \begin{equation}\label{G-Sylvester_over}
          AX-B^{\mathrm{T}}XS^{\mathrm{T}}=C-(SC)^{\mathrm{T}}.
        \end{equation}
      \end{jthrm}
      \begin{jprf}
        First, we show that $G:=I_m\otimes I_m-P_{mm}(I_m\otimes S)$ is nonsingular.
        Then we prove that the T-congruence Sylvester equation~\eqref{T-Sylvester} and the generalized Sylvester equation~\eqref{G-Sylvester_over} are equivalent.\par
        Let $(\lambda^{(G)}, \bm{v})$ be an eigenpair of $G$.
        Then it follows that
        \begin{align}
          P_{mm}(I_m\otimes S)\bm{v}=(1-\lambda^{(G)})\bm{v}.
        \end{align}
        Thus, $G$ is nonsingular if and only if $K:=P_{mm}(I_m\otimes S)$ has no eigenvalue equal to 1.
        From~\eqref{permu_1} and~\eqref{permu_4}, it follows that
        \begin{align}
          K^2&=P_{mm}(I_m\otimes S)P_{mm}(I_m\otimes S)\\
          &=(S\otimes I_m)(I_m\otimes S)\\
          &=S\otimes S.\label{SS}
        \end{align}
        As a result of~\eqref{SS} and Lemma~\ref{kronecker_eigenvalue}, any eigenvalue of $K^2$ is a product of two eigenvalues of $S$.
        Therefore, thanks to the assumption that $\Lambda(S)$ is reciprocal free, $G$ is nonsingular.\par
        Now, we show the statement of Theorem~\ref{overdetermined}.
        By multiplying~\eqref{linear_eq} by $G$ and using~\eqref{permu_1} and~\eqref{permu_4}, we have
        \begin{align}
          \{I_m\otimes A-P_{mm}(I_m\otimes SA)+P_{mm}(I_m\otimes B^{\mathrm{T}})-P_{mm}(I_m\otimes S)P_{mm}(I_m\otimes B^{\mathrm{T}})\}\bm{x}=G\bm{c},
        \end{align}
        which implies
        \begin{align}
          (I_m\otimes A-S\otimes B^{\mathrm{T}})\bm{x}=G\bm{c},\label{vec_overGseq}
        \end{align}
        since $B^\mathrm{T}=SA$.
        From~\eqref{permu_3} and~\eqref{permu_4}, the right-hand side of~\eqref{vec_overGseq} becomes
        \begin{align}
          G\bm{c}&=\bm{c}-P_{mm}(I_m\otimes S)\bm{c}\\
          &=\bm{c}-(S\otimes I_m)P_{mm}\bm{c}\\
          &=\mathrm{vec}(C)-(S\otimes I_m)\mathrm{vec}(C^{\mathrm{T}}).
        \end{align}
        Therefore, applying the inverse vec operator to~\eqref{vec_overGseq}, we obtain~\eqref{G-Sylvester_over} by using Lemma~\ref{vec_ABC}.
        Hence, we complete the proof of Theorem~\ref{overdetermined}.
        \qed
      \end{jprf}
      In Theorem~\ref{overdetermined}, we have not mentioned the existence of the matrix $S$ satisfying the condition $B^{\mathrm{T}}=SA$.
      If the matrix $A$ has full column rank, there always exist infinitely many matrices $S$ satisfying $B^{\mathrm{T}}=SA$.
      The simplest one is expressed by using the Moore--Penrose inverse of $A$.
      \begin{jpro}\label{pro_over} 
        Let $m\geq n$, $A\in\mathbb{R}^{m\times n}$ have full column rank and $B\in\mathbb{R}^{n\times m}$.
        Then $S:=B^{\mathrm{T}}A^\dagger$ yields $B^{\mathrm{T}}=SA$.
      \end{jpro}
      \begin{jprf}
        When $A$ has full column rank, the Moore--Penrose inverse of $A$ can be written as
        \begin{equation}
          A^{\dagger}=(A^{\mathrm{T}}A)^{-1}A^{\mathrm{T}},
        \end{equation}
        which implies $B^{\mathrm{T}}=SA$.
      \end{jprf}
      Notice that  if $m=n$, Theorem \ref{overdetermined} with Proposition~\ref{pro_over}  is equivalent to Theorem~\ref{moth}. Thus it can be regarded as an extension of  Theorem~\ref{moth}.
      Theorem~\ref{overdetermined} also holds true for the case $m<n$; 
      however,   $B^{\mathrm{T}}=SA$ holds true only for limited cases, i.e., an additional assumption is  required for $B$ and $A$.
      Thus, in the next subsection, we consider another transformation.\par 

      \subsection{The case $m\leq n$} 
      In this subsection, we describe a  novel transformation from the T-congruence Sylvester equation~\eqref{T-Sylvester} to the generalized Sylvester equation in the case $m\leq n$.
      \begin{jthrm}\label{underdetermined}
        Let $m\leq n$, $A\in\mathbb{R}^{m\times n}$, $B\in\mathbb{R}^{n\times m}$, and $C\in\mathbb{R}^{m\times m}$.
        Assume that there exists a matrix $D\in\mathbb{R}^{n\times m}$ such that $I_m=AD$, and for $S:=B^{\mathrm{T}}D$, $\Lambda(S)$ is reciprocal free.
        Then the T-congruence Sylvester equation~\eqref{T-Sylvester} is equivalent to the generalized Sylvester equation
        \begin{equation}\label{G-Sylvester_eq_under}
          A\tilde{X}-B^{\mathrm{T}}\tilde{X}S^{\mathrm{T}}=C,
        \end{equation}
        where $\tilde{X}\in\mathbb{R}^{n\times m}$ satisfies $X=\tilde{X}-D\tilde{X}^{\mathrm{T}}B$.
      \end{jthrm}
      \begin{jprf}
        First, we show that $G:=I_m\otimes I_n-P_{nm}(D\otimes B^{\mathrm{T}})$ is nonsingular.
        Then we prove that the T-congruence Sylvester equation~\eqref{T-Sylvester} and the generalized Sylvester equation~\eqref{G-Sylvester_eq_under} are equivalent.\par
        Let $(\lambda^{(G)}, \bm{v})$ be an eigenpair of $G$.
        Then, it follows that
        \begin{align}
          P_{nm}(D\otimes B^{\mathrm{T}})\bm{v}=(1-\lambda^{(G)})\bm{v}.
        \end{align}
        Thus, $G$ is nonsingular if and only if $K:=P_{nm}(D\otimes B^{\mathrm{T}})$ has no eigenvalue equal to 1.
        From~\eqref{permu_1} and~\eqref{permu_4}, it follows that
        \begin{align}
          K^2&=P_{nm}(D\otimes B^{\mathrm{T}})P_{nm}(D\otimes B^{\mathrm{T}})\\
          &=(B^{\mathrm{T}}\otimes D)(D\otimes B^{\mathrm{T}})\\
          &=B^{\mathrm{T}}D\otimes DB^{\mathrm{T}}. \label{BDDB}
        \end{align}
        Since $B^{\mathrm{T}}D$ and $DB^{\mathrm{T}}$ have the same set of nonzero eigenvalues, then we have from~\eqref{BDDB} and Lemma~\ref{kronecker_eigenvalue} that a square of any eigenvalue of $K$ is 0 or a product of two eigenvalues of $S$.
        Therefore, from the assumption that $\Lambda(S)$ is reciprocal free, $G$ is nonsingular.\par
        Now, we prove the statement of Theorem~\ref{underdetermined}.
        By using the nonsingular matrix $G$,~\eqref{permu_1},~\eqref{permu_2}, and~\eqref{permu_4}, the equation~\eqref{linear_eq} becomes
        \begin{align}
          \{I_m\otimes A+P_{mm}(I_m\otimes B^{\mathrm{T}})\}\{I_m\otimes I_n-P_{nm}(D\otimes B^{\mathrm{T}})\}\tilde{\bm{x}}=\bm{c},
        \end{align}
        which implies
        \begin{align}
          (I_m\otimes A-S\otimes B^{\mathrm{T}})\tilde{\bm{x}}=\bm{c}, \label{vec_underGseq}
        \end{align}
        where $\tilde{\bm{x}}:=G^{-1}\bm{x}$. 
        Since $\mathrm{vec}(\tilde{X})=\tilde{\bm{x}}$, applying the inverse vec operator to~\eqref{vec_underGseq}, we obtain~\eqref{G-Sylvester_eq_under} by using Lemma~\ref{vec_ABC}.
        Hence we complete the proof of Theorem~\ref{underdetermined}.\qed
      \end{jprf}
      In Theorem~\ref{underdetermined}, there exist infinitely many matrices $D$.
      The simplest one is the Moore--Penrose inverse of $A$.
      \begin{jpro}\label{pro_under}
        Let $m\leq n$, $A\in\mathbb{R}^{m\times n}$ have full row rank and $B\in\mathbb{R}^{n\times m}$.
        Then $D:=A^\dagger$ yields $I_m=AD$.
      \end{jpro}
      \begin{jprf}
        When $A$ has full row rank, the Moore--Penrose inverse of $A$ can be written as
        \begin{equation}
          A^{\dagger}=A^{\mathrm{T}}(AA^{\mathrm{T}})^{-1}.
        \end{equation}
        Thus, letting $D=A^{\dagger}$, we have $AD=I_m$.
      \end{jprf}
      In the square case, by setting $\hat{X}=A\tilde{X}$, Theorem~\ref{underdetermined} reduces to the following result:
      \begin{jcor}\label{cor_under}
        Let $A, B, C\in\mathbb{R}^{n\times n}$.
        Assume that $A$ is nonsingular, and for $S:=B^{\mathrm{T}}A^{-1}$, $\Lambda(S)$ is reciprocal free.
        Then, the T-congruence Sylvester equation~\eqref{T-Sylvester} is equivalent to the Lyapunov equation
        \begin{equation}\label{Lyapunov_eq_under}
          \hat{X}-S\hat{X}S^{\mathrm{T}}=C,
        \end{equation}
        where $\hat{X}:=A\tilde{X}$ and $\tilde{X}$ satisfies $X=\tilde{X}-A^{-1}\tilde{X}^{\mathrm{T}}B$.
      \end{jcor}
      This is a different result from that given in Theorem~\ref{moth}.
    \section{Conclusion}
      Our results are summarized in Table~\ref{table}. The key result is that under certain conditions T-congruence Sylvester equation is equivalent to  the generalized Sylvester equation.
      This may provide flexibility in research for finding further mathematical features and solvers of T-congruence Sylvester equation (or the generalized Sylvester equation) because mathematical features and solvers of the generalized Sylvester equation can be utilized for T-congruence Sylvester equation, and vice versa.
      \par
      \begin{table}[hbtp]
        \caption{Equivalent matrix equations of the T-congruence Sylvester equation $AX+X^{\mathrm{T}}B=C.$}\label{table}
        \centering
        \begin{tabular}{c||c|c}
          \hline
          &$m\geq n$&$m\leq n$\\
          \hline\hline
          Theorem
          &\begin{tabular}{c}
            generalized Sylvester eq. (Th. 3.2)\\
            $AX-B^{\mathrm{T}}XS^{\mathrm{T}}=C-(SC)^{\mathrm{T}}$
          \end{tabular}
          &\begin{tabular}{c}
            generalized Sylvester eq. (Th. 3.5)\\
            $A\tilde{X}-B^{\mathrm{T}}\tilde{X}S^{\mathrm{T}}=C$
          \end{tabular}\\
          \hline
          $m=n$&\begin{tabular}{c}
            Lyapunov eq.~\cite[Th. 3.1]{oozawa2018}\\
            $\tilde{X}-S\tilde{X}S^{\mathrm{T}}=Q$
          \end{tabular}
          &\begin{tabular}{c}
            Lyapunov eq. (Cor.~\ref{cor_under})\\
            $\hat{X}-S\hat{X}S^{\mathrm{T}}=C$
          \end{tabular}\\
          \hline
          \end{tabular}
        \end{table}





\end{document}